\documentclass[11pt]{article}

\usepackage{amsmath}
\usepackage{amsfonts}
\usepackage{hyperref}
\usepackage{graphicx}
\usepackage{multirow}
\usepackage{xcolor}
\usepackage{latexsym}        
\usepackage{enumitem}
\usepackage{ushort}
\usepackage{accents}
\usepackage{soul}

\newtheorem{teo}{Theorem}[section]
\newtheorem{lem}{Lemma}[section]
\newtheorem{coro}{Corollary}[section]
\newtheorem{assump}{Assumption}[section]
\newcommand{\halmos}{\hfill$\Box$}
\newenvironment{pro}{\noindent\textit{Proof:}}{\halmos}
\newtheorem{defi}{Definition}[section]

\newcommand{\xtrial}{x_\mathrm{trial}}

\newcommand{\R}{\mathbb{R}}
\newcommand{\N}{\mathbb{N}}
\newcommand{\cm}{\mathrm{cm}}

\newcommand{\yrest}{ y^{\mathrm{re}}_k}

\newcommand{\hmax}{h_{\max}}
\newcommand{\fmin}{f_{\min}}

\newcommand{\feas}{\mathrm{feas}}
\newcommand{\opt}{\mathrm{dist}}

\newcommand{\alphabis}{\gamma}

\begin{document}

\title{Inexact restoration for derivative-free expensive
  function minimization and applications \thanks{This work has been
    partially supported by the Brazilian agencies FAPESP (grants
    2013/07375-0, 2016/01860-1, and 2018/24293-0) and CNPq (grants
    302538/2019-4 and 302682/2019-8) and by the Serbian Ministry of
    Education, Science, and Technological Development.}}

\author{
  E. G. Birgin\thanks{Department of Computer Science, Institute of
    Mathematics and Statistics, University of S\~ao Paulo, Rua do
    Mat\~ao, 1010, Cidade Universit\'aria, 05508-090, S\~ao Paulo, SP,
    Brazil. e-mail: egbirgin@ime.usp.br}
  \and
  N. Kreji\'c\thanks{Department of Mathematics and Informatics,
    Faculty of Sciences, University of Novi Sad, Trg Dositeja
    Obradovi\'ca 4, 21000 Novi Sad, Serbia. e-mail:
    natasak@uns.ac.rs}
  \and
  J. M. Mart\'{\i}nez\thanks{Department of Applied Mathematics,
    Institute of Mathematics, Statistics, and Scientific Computing
    (IMECC), University of Campinas, 13083-859 Campinas SP,
    Brazil. e-mail: martinez@ime.unicamp.br}}

\date{June 3, 2021}

\maketitle

\begin{abstract}
The Inexact Restoration approach has proved to be an adequate tool for
handling the problem of minimizing an expensive function within an
arbitrary feasible set by using different degrees of precision. This
framework allows one to obtain suitable convergence and complexity
results for an approach that rationally combines low- and
high-precision evaluations.  In this paper we consider the case where
the domain of the optimization problem is an abstract metric
space. Assumptions about differentiability or even continuity will not
be used in the general algorithm based on Inexact Restoration.
Although optimization phases that rely on smoothness cannot be used in
this case, basic convergence and complexity results are recovered. A
new derivative-free optimization phase is defined and the subproblems
that arise at this phase are solved using a regularization approach
that takes advantage of different notions of stationarity.  The new
methodology is applied to the problem of reproducing a controlled
experiment that mimics the failure of a dam.\\

\noindent
\textbf{Key words:} Nonlinear programming, Inexact Restoration,
derivative-free, inexact evaluation of expensive function,
algorithms.\\

\noindent
\textbf{2010 Mathematics Subject Classification:} 65K05, 65K10,
90C30, 90C56, 65Y20, 90C90.
\end{abstract}

\section{Introduction}

For many reasons scientists and engineers may need to optimize
problems in which the objective function is very expensive to
evaluate. In these cases, partial, and obviously inexact, evaluations
are useful. The idea is to decrease as much as possible functional
values using partial evaluations in such a way that, when we have no
choice except to evaluate the function with maximal accuracy, we are
already close enough to a solution of the problem. Rational decisions
about when to increase accuracy (and evaluation cost) or even when to
try more inexact evaluations are hard to make. Roughly speaking, we
need a compromise between accuracy of evaluation and functional
decrease that is difficult to achieve on a mere heuristic basis.

Most papers on minimization methods with inexact evaluations aim to
report the behavior of modifications of standard algorithms in the
presence of errors in the computation of the objective function, its
derivatives, or the constraints
\cite{inexact1,inexact2,inexact3,inexact4,inexact5,inexact6,inexact7}.
In general, it is assumed that the objective function and, perhaps
its derivatives, can be computed with a given error bound each time
the functional value is required at an arbitrary
point. In~\cite{inexact1}, the underlying ``exact'' algorithm is based
on adaptive regularization. Trust-region algorithms inspired
in~\cite[\S10.6]{cgtbook} are addressed in~\cite{inexact6}. High-order
complexity of stochastic regularization is considered
in~\cite{inexact2}.  The case of PDE-constrained optimization is
studied in \cite{inexact3}.  Nonsmoothness is addressed
in~\cite{inexact4}. Convex quadratic problems using Krylov methods are
studied in \cite{inexact7}. Many of these works seem to be influenced
by an early paper by Carter \cite{inexact5}.

In recent papers~\cite{bkm2018,bkm2019,nkjmm}, a methodology based on
the analogy of the Inexact Restoration idea for continuous constrained
optimization and the process of increasing accuracy of function
evaluations was developed. In this approach, it is not assumed that
one is able to compute a functional value (let alone derivatives)
within a given required error bound. Instead, it is assumed that
accuracy is represented by an abstract function $h(y)$ such that
$h(y)=0$ means maximal accuracy and $y$ is a case-dependent procedure.
For instance, $y \in Y$ may represent an algorithm with which one can
compute the approximate objective function and $h(y) \geq 0$ is an
accuracy-related function. The connection between the value of $h(y)$
and a possible error bound is not assumed to be known. For example,
$y$ may represent the maximal number of iterations that are allowed
for a numerical algorithm that computes the objective function before
obtaining convergence when we know that convergence eventually occurs
but an error estimation is not available. Several additional examples
may be found in \cite{bkm2019}.
%It has been verified that Inexact
%Restoration provides an adequate framework to deal with this case.

Inexact Restoration methods for smooth constrained optimization were
introduced in~\cite{ir2}. Each iteration of an Inexact Restoration
method proceeds in two stages. In the Restoration Phase, infeasibility
is reduced; and, in the Optimization Phase, the reduction of the
objective function or its Lagrangian is addressed with a possible loss
of feasibility. Convergence with sharp Lagrangians as merit functions
was proved in~\cite{ir4}. Applications to bilevel programming were
given in~\cite{ir5,ir7}. In \cite{ir19}, an application to
multiobjective optimization was described. The first line-search
implementation was introduced in~\cite{ir3}. Nonsmooth versions of the
main algorithms were defined in~\cite{ir1,ir17,ir18}. The employment
of filters associated with Inexact Restoration was exploited
in~\cite{ir6,ir17}. Applications to control problems were given
in~\cite{ir8,ir9,ir10,ir13}. In~\cite{ir11}, Inexact Restoration was
used to obtain global convergence of a sequential programming
method. Large-scale applications were discussed
in~\cite{ir12}. In~\cite{ir14}, Inexact Restoration was used for
electronic structure calculations; and problems with a similar
mathematical structure were addressed in~\cite{ir27}. The reliability
of Inexact Restoration for arbitrary nonlinear optimization problems
was assessed in \cite{ir15}. In \cite{ir22}, the worst-case evaluation
complexity of Inexact Restoration was analyzed. An application to
finite-sum minimization was described in~\cite{ir23}. Continuous and
discrete variables were considered in~\cite{ir24}. Non-monotone
alternatives were defined in~\cite{ir25}. An application to the demand
adjustment problem was given in~\cite{ir26}. Local convergence results
were proved in~\cite{ir9}.

In~\cite{bkm2019}, an algorithm of Inexact Restoration type applied to
minimization with inexact evaluations that exhibits convergence and
complexity results, was developed. However, the algorithm proposed
in~\cite{bkm2019}, as well as the algorithms previously introduced
in~\cite{bkm2018,nkjmm}, employs derivatives of the objective
function, a feature that may be inadequate in many cases in which one
does not have differentiability at all. This state of facts motivates
the present work. Here, the algorithms of~\cite{bkm2018,bkm2019,nkjmm}
are adapted to the case in which derivatives are not available and the
main theoretical results are proved. In addition, the domain of the
optimization problem is considered to be a metric space, unlike an
Euclidean space as considered in~\cite{bkm2018,bkm2019,nkjmm}. The
reason why we are proposing the generalization to metric spaces is
simple: Many real optimization problems take place in domains that are
not Euclidean spaces and, many times, in domains that are not even
vector spaces. Naturally, there is a price we pay for this
generalization because there are no known optimality conditions that
can be invoked to define possible approximate solutions to the
problems. In the current state of development of methods for
optimization in metric spaces, it is natural that the optimality
conditions that we can define are inspired by what happens in
Euclidean spaces. In that sense, we emphasize a poorly known or poorly
exploited property: In the Euclidean case, a point being a local
minimizer does not imply that that point is also a local minimizer of
a high-order approximation of the objective function. (For example,
$(0, 0)$ is a minimizer of $x_2^2 - x_1^3 x_2 + x_1^6$ but is not a
local minimizer of its Taylor approximation $x_2^2 - x_1^3 x_2$.)
That implication is only valid when the dimension of the space is 1 or
when the approximation is of order 2. If we deal with higher order
approximations, the fact that a point is a local minimizer only
implies that it is also a local minimizer of a suitably regularized
Taylor approximation.  That is the reason why, in the case of metric
spaces, we cannot establish stronger necessary conditions than those
mentioned in the present article. With these limitations, we elaborate
the possible theory that frames affordable and reasonably efficient
techniques. Moreover, in the present work, we focus on a practical
problem related with the prediction and mitigation of the consequences
of dam breaking disasters.

The rest of this paper is organized as follows. The main contributions
of the paper are surveyed in Section~\ref{contributions}.  The basic
Inexact Restoration algorithm for problems with inexact evaluation
without error bounds is introduced in Section~\ref{mainalgo}. This
algorithm is the one described in~\cite{bkm2019} with the difference
that the domain $\Omega$ is an arbitrary metric space here, instead of
a subset of $\R^n$. This extension may be useful for cases in which we
have discrete or qualitative variables or when the original problem is
formulated on a functional space. The proofs of the complexity and
convergence results for this algorithm are similar to the ones
reported in~\cite{bkm2019}.  Since the employment of derivatives is
essential for the definition of the Optimization Phase
in~\cite{bkm2019}, in the present contribution we need a different
approach for that purpose.  This is the subject of
Section~\ref{solvinsub}. The Optimization Phase requires sufficient
descent of an iteration-dependent inexactly evaluated
function. Sufficient descent of this function is equivalent to simple
descent of an associated forcing function, therefore the required
sufficient descent may be obtained by applying an arbitrary number of
iterations of an \textit{ad hoc} monotone minimization algorithm to
an associated forcing function. This monotone algorithm, that depends
on characteristics of the function being minimized, is called MOA
(Monotone Optimization Algorithm) in the rest of the paper. In this
phase, we may take advantage of the theoretically justified
possibility of relaxing precision.  In Section~\ref{regula}, we
concentrate at the definition of a strategy for MOA. Our choice is to
solve the problem addressed by MOA by means of suitable
surrogate models associated with iterated adaptive
regularization. Theoretical results of the previous sections are
summarized in Section~\ref{condensed}. In Section~\ref{experiments},
we describe an implementation that is adequate for the practical
problem that we have in mind. The idea is to simulate a small-scale
controlled physical experiment that mimics the failure of a dam,
described in \cite{pirulli}. We rely on an MPM-like (Moving Particles
Method) approach in which the dynamic of the particles is governed by
the behavior of the Spectral Projected Gradient (SPG) method applied
to the minimization of a semi-physical energy. Our optimization
problem consists of finding the parameters of the energy function by
means of which the trajectory of the SPG method best reproduces
physical experiments reported in~\cite{pirulli}.  In
Section~\ref{conclusions}, we state conclusions and lines for future
research.\\

\noindent
\textbf{Notation.} $\N_+$ denotes the non-negative integer numbers;
while $\R_+$ denotes the non-negative real numbers. The symbol $\|
\cdot \|$ denotes the Euclidean norm of vector and matrices.

\section{Main problem and contributions} \label{contributions}

Consider the problem
\begin{equation} \label{theproblem}
\mbox{Minimize (with respect to $x$) } f(x,y)
\mbox{ subject to } h(y)=0 \mbox{ and } x \in \Omega,
\end{equation}
where $\Omega$ is an arbitrary metric space, $Y$ is an abstract set,
$h: Y \to \R_+$, and $f: \Omega \times Y \to \R$. For further
reference, $d(x,z)$ denotes the distance between $x$ and $z$ in
$\Omega$. $Y$ is sometimes interpreted as a set of indices according
to which the objective function is computed with different levels of
accuracy. The nature of these indices varies with the problem and, for
this reason, $Y$ is presented as an arbitrary set.  The function
$h(y)$ is assumed to be related with the accuracy in the evaluation of
$f$ associated with the index $y \in Y$. The smaller the value of
$h(y)$, the higher the accuracy; while maximal accuracy corresponds to
$h(y) = 0$. In principle, problem~(\ref{theproblem}) consists of
minimizing the objective function with maximal accuracy.

In this article, the domain~$\Omega$ of the optimization
problem~(\ref{theproblem}) is a metric space, opposed to our previous
paper~\cite{bkm2019}, where it was the Euclidean space or a subset of
it. We have in mind, therefore, spaces formed by functions in which it
is possible to define a metric, and optimization problems defined on
discrete sets. Because of this level of generality, the objective
function~$f$ can rarely be computed exactly, so inaccurate evaluations
are unavoidable.  However, there are many ways to deal with inaccuracy
of functional evaluations.  In this article each instance of
inaccurate evaluation is identified by an element $y$ of an abstract
set $Y$. For example, if the evaluation of the function comes from an
iterative process, each $y$ can be identified with the maximum number
of iterations allowed. If the evaluation depends on a discretization
grid, $Y$ is the set of all possible grids.  To each $y$, a ``score''
$h(y)$ is attributed such that $h(y)=0$ corresponds to the exact
evaluation of the objective function while $h(y) > 0$ corresponds to
different degrees of inaccuracy; and the smaller the value of $h(y)$,
the higher the evaluation accuracy.  Therefore,
problem~(\ref{theproblem}) consists in minimizing $f(x, y)$ with
respect to $x \in \Omega$, restricted to $h(y)=0$.

The evaluation of $f(x,y)$ when $h(y)=0$ is very often impossible or
extremely expensive, therefore, we must rely on evaluations of
$f(x,y)$ with $h(y)$ small. However, the cost of evaluating $f(x,y)$
generally increases as $h(y)$ decreases. Therefore, the strategy for
solving problem~(\ref{theproblem}) relies on solving problems of the
form ``Minimize $f(x, y_k)$'' with $h(y_k) \to 0$. A rational way of
doing this has been defined in~\cite{bkm2019} for the case in which
the domain is Euclidean and is extended here to the general case.

Clearly, with this level of generality differentiability arguments are
out of question and we need ``derivative-free'' procedures for the
(approximate) minimization of $f(x,y_k)$.  For this purpose we
consider that the objective function may be approximated, at each
iteration, by a suitable model and that approximate minimization of
such model is affordable.  The precise determination of the model is,
of course, problem-dependent. Under suitable assumptions, we prove
that, in finite time, it is possible to obtain good approximations to
a solution of problem~(\ref{theproblem}). Finally, we present a
practical application that concerns fitting models of dam breaking.

\section{Main algorithm} \label{mainalgo}

In this section, the problem, the algorithm, and the convergence
results introduced in~\cite{bkm2019} are revisited. Of note is that in
the present contribution the domain of the problem is a general metric
space instead of a subset of $\R^n$. The proofs of the theorems
presented in this section are analogous to the ones given in
\cite{bkm2019}, with the occurrences of $\|x - z\|$ replaced by the
distance $d(x,z)$ between $x$ and $z$.

Let us define a merit function $\Phi: \Omega
\times Y \times (0, 1) \to \R$ by
\[
\Phi(x, y, \theta) = \theta f(x, y) + (1-\theta) h( y ).
\]   
The merit function $\Phi$ combines objective function evaluation and
accuracy level. Here, $f(x,y)$ is the value of the objective function
for a feasible $x$ when the accuracy level is determined by $y \in
Y$. Therefore, $\Phi(x,y,\theta)$ represents a compromise between
optimality and accuracy. The penalty parameter $\theta$ is updated at
each iteration of the main algorithm. The main algorithm, with an
\textit{Optimization Phase} different from the one of the algorithm
introduced in~\cite{bkm2019}, follows.\\

\noindent
\textbf{Algorithm~\ref{mainalgo}.1}. Let $x_0 \in \Omega$, $y_0 \in Y$,
$\theta_0 \in (0,1)$, $\nu > 0$, $r \in (0,1)$, $\alpha > 0$, and
$\beta > 0$ be given. Set $k \leftarrow 0$.
\begin{description}
\item[\textbf{Step~1.}] \textit{Restoration phase}

Define $\yrest \in Y$ in such a way that
\begin{equation} \label{erre}
h(\yrest) \leq r h(y_k)
\end{equation}
and
\begin{equation} \label{beta}
f(x_k, \yrest) \leq f(x_k, y_k) + \beta h(y_k).
\end{equation}

\item[\textbf{Step~2.}] \textit{Updating the penalty parameter}

If
\begin{equation} \label{cinco}
\Phi(x_k, \yrest, \theta_k) \leq \Phi(x_k, y_k, \theta_k) +
\frac{1-r}{2} \left( h(\yrest) - h(y_k) \right),
\end{equation}
set $\theta_{k+1} = \theta_k$. Otherwise, set
\begin{equation} \label{seis}
\theta_{k+1} = \frac {(1+r) \left( h(y_k) - h(\yrest) \right)} {2
  \left( f(x_k, \yrest) - f(x_k, y_k) + h(y_k) - h( \yrest ) \right)}.
\end{equation}

\item[\textbf{Step~3.}] \textit{Optimization phase}

Compute $ y_{k+1} \in Y$ and $x_{k+1} \in \Omega$, such that 
\begin{equation} \label{armijo1}
f(x_{k+1}, y_{k+1}) \leq f(x_k, \yrest) - \alpha d(x_k, x_{k+1})^\nu
\end{equation}     
and
\begin{equation} \label{dos}
\Phi(x_{k+1}, y_{k+1}, \theta_{k+1}) \leq \Phi(x_k, y_k,
\theta_{k+1}) + \frac{1-r}{2} \left( h(\yrest) - h(y_k) \right).
\end{equation} 
Update $k \leftarrow k+1$, and go to Step~1.
\end{description}

Given an iterate $(x_k, y_k)$, at Step~1 we compute the function with
an accuracy determined by $\yrest \in Y$, which is better than the
accuracy given by $y_k$. At Step~2 we update the penalty parameter
with the aim of decreasing the merit function at the restored pair
$(x_k, \yrest)$ with respect to the merit function computed at $(x_k,
y_k)$. This is typical in Inexact Restoration methods for constrained
optimization. Note that, using (\ref{erre}), (\ref{beta}), the
negation of~(\ref{cinco}), and the fact that $\theta_k > 0$, it is
easy to show that both the numerator and the denominator
of~(\ref{seis}) are positive and that $\theta_{k+1} \in (0,
\theta_k]$. At Step~3 (optimization phase), we compute the new iterate,
at which the objective function value should decrease with respect to
the restored pair in the sense of (\ref{armijo1}) and the merit
function should be improved with respect to $(x_k, y_k)$ in the sense
of (\ref{dos}).

%The main results in this section are generalizations of the ones
%proved in \cite{bkm2019}. 

The theoretical results below describe the properties of
Algorithm~\ref{mainalgo}.1 and the generated sequence of
iterates. These results state the algorithm is well-defined and that
any desired accuracy level, for the evaluation of the objective
function, can be reached in a finite number of iterations. (In fact, a
complexity result is presented and an upper bound on the number of
iterations is given.) A first step towards an optimality result,
including complexity bounds, is also given. The optimality results
will be completed in the remaining of the work with the definition of
a particular strategy for computing~$x_{k+1}$ and~$y_{k+1}$ at Step~3
and the analysis of its convergence properties.

In Assumption~\ref{a1} we state that defining~$\yrest$
satisfying~(\ref{erre}) and~(\ref{beta}) at Step~1 of
Algorithm~\ref{mainalgo}.1 is always possible. Notice that assuming we
can define~$\yrest$ satisfying~(\ref{erre}) is a simple assumption,
since it corresponds to assuming that we have the capability of
increasing the accuracy for the evaluation of the objective
function. On the other hand, the capability of making the
defined~$\yrest$ to satisfy~(\ref{beta}) is problem-dependent; and
plausible choices are given in~\cite{bkm2019}. In Assumption~\ref{a2},
basic assumptions on functions~$h$ and $f$ are stated.

\begin{assump} \label{a1}
At Step~1 of Algorithm~\ref{mainalgo}.1, for all $k \in \N_+$ it is
possible to compute, in finite time, $\yrest$ satisfying~(\ref{erre})
and (\ref{beta}).
\end{assump}

\begin{lem} \label{lemaboludo}
Suppose that Assumption~\ref{a1} holds. Then,
Algorithm~\ref{mainalgo}.1 is well defined.
\end{lem}    

\begin{pro}
Replace, in Lemma~2.3 of \cite{bkm2019}, the occurrences of $\|x-z\|$
with $d(x, z)$.
\end{pro}

\begin{assump} \label{a2}
There exist $\hmax > 0$ and $\fmin \in \R$ such that, for all $y \in
Y$ and $x \in \Omega$ we have that $h(y) \leq h_{\max}$ and $f(x, y)
\geq f_{\min}$.
\end{assump}
 
\begin{teo} \label{corofeas}
Suppose that Assumptions~\ref{a1} and~\ref{a2} hold. Given
$\varepsilon_{\feas} > 0$, the number of indices~$k$ such that $h(y_k)
> \varepsilon_{\feas}$ is bounded above by
\begin{equation} \label{lavarlosdientes}
\frac{c_{\feas}}{\varepsilon_{\feas}},
\end{equation}
where $c_{\feas}$ only depends on $x_0$, $y_0$, $r$, $\theta_0$,
$\beta$, $h_{\max}$, and $f_{\min}$.
\end{teo}

\begin{pro}
This proof follows as the one of Corollary 2.2 of \cite{bkm2019} with
the replacements stated in the proof of Lemma~\ref{lemaboludo}.
\end{pro}

\begin{teo} \label{corosk}
Suppose that Assumptions~\ref{a1} and~\ref{a2} hold. Then, the series
$\sum_{k=0}^\infty d(x_k, x_{k+1})^\nu$ is convergent. Moreover, given
$\varepsilon_{\opt} > 0$, the number of iterates $k$ at which $d(x_k,
x_{k+1}) > \varepsilon_{\opt}$ is not bigger than
\begin{equation} \label{buenoeso}
\frac{c_{\opt}}{\varepsilon_{\opt}^{\nu}}, 
\end{equation}
where $c_{\opt}$ only depends on $\alpha$, $x_0$, $y_0$, $r$,
$\theta_0$, $\beta$, $h_{\max}$, and $f_{\min}$.
\end{teo}

\begin{pro}
See Theorem 2.2 and Corollary 2.3 of \cite{bkm2019} and perform the
substitutions of the norm of differences with appropriate distances.
\end{pro}

\begin{teo} \label{teoneo}
Suppose that Assumptions~\ref{a1} and~\ref{a2} hold. Then,
\begin{equation} \label{limhy}
\lim_{k \to \infty} h(y_k) = 0.
\end{equation}
Moreover, if $\Omega$ is complete and $\nu \leq 1$, the sequence
$\{x_k\}$ is convergent.
\end{teo}

\begin{pro} 
The property (\ref{limhy}) follows from Theorem~\ref{corofeas}. By the
convergence of the series $\sum_{k=0}^\infty d(x_k, x_{k+1})^\nu$ with
$\nu \leq 1$ we have that $\{x_k\}$ is a Cauchy sequence. Since
$\Omega$ is complete it turns out that $\{x_k\}$ is convergent.
\end{pro}\\

The convergence of the sequence $\{x_k\}$ when $\nu \leq 1$ had not
been considered in \cite{bkm2019}. In general, in nonconvex
optimization, it is proved that every convergent subsequence of the
sequence generated by an algorithm has optimality properties but
convergence of the whole sequence is rarely obtained. The validity of
the theorem does not necessarily lead us to use $\nu \leq 1$ in
numerical computations, since, on the one hand, the convergence to
zero of the distance between consecutive iterates is, in general, a
good indication of optimality and, on the other hand, values of $\nu$
greater than~1 do not prevent the total sequence from being
convergent.

\section{A general framework for the Optimization Phase} \label{solvinsub}

In this section, we consider the implementation of Step~3 of
Algorithm~\ref{mainalgo}.1.  Obviously, if, for any given precision
parameter~$y_{k+1}$, we apply any monotone optimization algorithm
(MOA) to the minimization of $f(x, y_{k+1}) + \alpha d(x_k, x)^\nu$
subject to $x \in \Omega$, starting with $x_k$ as initial
approximation, we find $ x \in \Omega$ such that
\[
f(x, y_{k+1}) \leq f(x_k, y_{k+1}) - \alpha d(x_k, x)^\nu.
\]
This follows from the fact that even $x = x_k$ satisfies such
inequality. However this is not a guarantee for the fulfillment of
\begin{equation} \label{moa1}
f(x, y_{k+1}) \leq f(x_k, \yrest) - \alpha d(x_k, x)^\nu.
\end{equation}  
and
\begin{equation} \label{moa2}
\Phi(x, y_{k+1}, \theta_{k+1}) \leq \Phi(x_k, y_k, \theta_{k+1}) +
\frac{1-r}{2} \left( h(\yrest) - h(y_k) \right),
\end{equation}
as required at Step~3 of Algorithm~\ref{mainalgo}.1. Therefore, before
accepting the ``solution'' obtained by MOA as new iterate $x_{k+1}$,
we must test whether (\ref{moa1}) and~(\ref{moa2}) hold. If
(\ref{moa1}) and~(\ref{moa2}) hold, then we accept $x_{k+1} =
x$. Otherwise, if after a finite number of trials with arbitrary
$y_{k+1}$ we fail to find $x$ satisfying (\ref{moa1})
and~(\ref{moa2}), we discard the MOA solution $x$ and apply MOA to the
minimization of $f(x, \yrest) + \alpha d(x_k, x)^\nu$. This time, the
``solution'' obtained by MOA necessarily yields the fulfillment of
(\ref{moa1}) and (\ref{moa2}) and we may accept $y_{k+1} = \yrest$ and
$x_{k+1}=x$. Two natural questions arise:
\begin{enumerate}
\item Why not using directly $y_{k+1} = \yrest$, instead of trying
  first one or more different precision parameters~$y_{k+1}$?
\item Which stopping criterion should we use for MOA?
\end{enumerate}
The answer to the first question is the following: We normally want to
choose $y_{k+1} \neq \yrest$ in such a way that $h(y_{k+1}) >
h(\yrest)$ because $h(y') > h(y)$ generally implies that the cost of
evaluating the objective function with precision given by $y'$ is
smaller than the cost of evaluating the objective function with
precision defined by $y$. In other words, we want to obtain the
maximum possible progress without unnecessary evaluation effort. The
answer to the second question is problem-dependent. Roughly speaking,
we assume that MOA is equipped with a stopping criterion related to
optimality. Recall that we are dealing with metric spaces in which
nothing similar to KKT conditions exists. The idea of using an
auxiliary algorithm for the fulfillment of descent conditions in the
Optimization Phase of Inexact Restoration comes from~\cite{bfms},
where Generating Set Search algorithms~\cite{klt,kltsiamreview} were
employed as subalgorithms for an Inexact Restoration method for
derivative-free optimization with smooth constraints.

Summing up, the strategy to compute $y_{k+1} \in Y$ and $x_{k+1} \in
\Omega$ at Step~3 of Algorithm~\ref{mainalgo}.1 is described by
Algorithm~\ref{solvinsub}.1 below.\\

\noindent
\textbf{Algorithm~\ref{solvinsub}.1.}
\begin{description}  
\item[\textbf{Step~1.}] Choose $y_{k+1} \in Y$.
\item[\textbf{Step~2.}] Define, for all $x \in \Omega$,
\begin{equation} \label{Efe}
F(x) = f(x, y_{k+1}) + \alpha d(x, x_k)^\nu.
\end{equation}

\item[\textbf{Step~3.}] Consider the subproblem
\begin{equation}\label{minefe}
\mbox{Minimize } F(x) \mbox{ subject to } x \in \Omega.
\end{equation}
By using an appropriate monotone iterative (derivative-free)
minimization Algorithm~MOA, starting with $x = x_k$, compute $\xtrial
\in \Omega$ such that
\begin{equation} \label{armijo8}
F(\xtrial) \leq F(x_k)
\end{equation}
and $\xtrial$ satisfies a stopping criterion related to MOA  
to be specified later.
\item[\textbf{Step 4.}] If $y_{k+1} = \yrest$ or 
\begin{equation} \label{67}
\begin{array}{c}
f(\xtrial, y_{k+1}) \leq f(x_k, \yrest) - \alpha d(x_k, \xtrial)^\nu \mbox{ and }\\[4mm]
\Phi(\xtrial, y_{k+1}, \theta_{k+1}) \leq \Phi(x_k, y_k,
\theta_{k+1}) + \frac{1-r}{2} \left( h(\yrest) - h(y_k) \right),
\end{array}
\end{equation}       
then \textbf{return} $y_{k+1}$ and $x_{k+1} = \xtrial $.
\item[\textbf{Step 5.}] Re-define $y_{k+1} = \yrest$, and go to
  Step~2.
\end{description}

\noindent
\textbf{Remark \ref{solvinsub}.1.} At Step~1 of
Algorithm~\ref{solvinsub}.1, we have the chance of choosing a looser
precision $y_{k+1}$, i.e.\ $y_{k+1}$ such that~$h(y_{k+1})>h(y_k)$,
and, with this precision, trying to improve (at Step~3) the current
approximate solution~$x_k$. The success of this attempt depends on the
capacity of computing $\xtrial$ satisfying~(\ref{armijo8}) at Step~3
that also satisfies~(\ref{67}). By means of a judicious choice of
$y_{k+1}$, we may obtain, simultaneously and with low computational
cost, a smaller function value and, consequently, substantial progress
in terms of distance to the true solution of the original problem.\\

\noindent
\textbf{Remark \ref{solvinsub}.2.} Note that (\ref{armijo8}) is
equivalent to
\begin{equation} \label{armijo2}
f(\xtrial, y_{k+1}) \leq f(x_k, y_{k+1}) - \alpha d(x_k, \xtrial)^\nu
\end{equation}    
and that when $y_{k+1} = \yrest$ the fulfillment of~(\ref{armijo2})
implies trivially the fulfillment of~(\ref{67}). This is the reason
why, at Step~4, the test of~(\ref{67}) is not necessary when $y_{k+1}
= \yrest$.\\

Returning back to Algorithm~\ref{mainalgo}.1, by
Theorems~\ref{corofeas} and~\ref{corosk}, given
$\varepsilon_{\feas}>0$ and $\varepsilon_{\opt}>0$, there exists an
iteration index~$k$ not larger than
\begin{equation} \label{existek}
\frac{c_{\feas}}{\varepsilon_{\feas}} +
\frac{c_{\opt}}{\varepsilon_{\opt}^{\nu}},
\end{equation}
where $c_{\feas}$ and $c_{\opt}$ are constants that only depend on the
problem and algorithmic constants, such that
\begin{equation} \label{hyd}
h(y_{k+1}) \leq \varepsilon_{\feas} \mbox{ and } d(x_{k}, x_{k +1})
\leq \varepsilon_{\opt}.
\end{equation}
It is quite natural to ask whether (\ref{hyd}) implies, in some sense,
feasibility and optimality of the original problem.  Observe, firstly,
that, by Step~3 (of Algorithms~\ref{mainalgo}.1
and~\ref{solvinsub}.1), $x_{k+1} \in \Omega$. Moreover, there is no
doubt that $h(y_{k+1}) \leq \varepsilon_{\feas}$ implies that the
objective function has been computed at $x_{k+1}$ with high precision.
Finally, $x_{k+1}$ is an approximate minimizer of $f(x, y_{k+1}) +
\alpha d(x, x_k)$ in the sense that this point satisfies an optimality
condition related to MOA. Since $\alpha$ (parameter of
Algorithm~\ref{mainalgo}.1) is assumed to be small and, by
(\ref{hyd}), $x_{k+1}$ is close to $x_k$, this is the best
justification that we have for saying that, in some sense, $x_{k+1}$
is almost optimal. Recall that we work in a general metric space where
a vectorial structure for defining optimality conditions is not
assumed to exist at all.  As a consequence, our confidence in the fact
that (\ref{existek}) and (\ref{hyd}) reflect optimality relies in the
properties of the algorithm MOA. This is the subject of the following
section.

\section{Monotone Optimization Algorithm} \label{regula}

At Step~3 of Algorithm~\ref{solvinsub}.1, we need to find $\xtrial \in
\Omega$ such that
\begin{equation} \label{armijo9}
F(\xtrial) \leq F(x_k)
\end{equation}
and $\xtrial$ satisfies a stopping criterion related to MOA, where MOA
is an appropriate monotone optimization algorithm applied to the
minimization of $F(x)$ subject to $x \in \Omega$. Unfortunately, there
are no implementable algorithms for minimizing a function over a
general metric space. Therefore, the problem of minimizing $F(x)$ for
$x \in \Omega$ must be solved by means of a reduction to a finite
dimensional linear space.  We perform this reduction in two steps.
Let us define, for all $\bar x, x \in \Omega$, a surrogate model
$M(\bar x, x)$ such that the following assumption holds.

\begin{assump} \label{a3}
For each $\bar x \in \Omega$, $M(\bar x, x)$ satisfies:
\begin{enumerate}
\item For all $\bar x \in \Omega$, $M(\bar x, \bar x) = F(\bar x)$.
\item There exists $L \geq 0$ and $ p + 1> 0$ such that, for all $\bar
  x, x \in \Omega$,
  \begin{equation} \label{lipscho}
    F(x) \leq M(\bar x, x) + L \, d(\bar x, x)^{p+1}.
  \end{equation}
\item There exists $\sigma \geq 0$ such that $M(\bar x, x) + \sigma
  d(\bar x, x)^{p+1}$ is bounded below in $\Omega$.
\item For all $\bar x \in \Omega$ and $\sigma > 0$, it is affordable
  to solve, up to some given precision, the problem
  \begin{equation} \label{surro} 
    \mbox{Minimize } M(\bar x, x) + \sigma d(\bar x, x)^{p+1}
    \mbox{ subject to } x \in \Omega.
  \end{equation}
\end{enumerate}      
\end{assump}

The idea is to substitute the (approximate) minimization of~$F(x)$
with the (approximate) minimization of~(\ref{surro}). It could be
argued that~(\ref{surro}) is a minimization problem over a metric
space, as well as the minimization of~$F(x)$. However, it is expected
that a model of~$F(x)$ can be built whose minimization over the metric
space reduces to a finite dimensional minimization. For example,
if~$\Omega$ is a space of continuous functions, we may approximate the
minimization over~$\Omega$ by a minimization over a finite dimensional
space of splines and we may employ Taylor-like models for this
minimization. As a whole, we could obtain the structure required by
Assumption~\ref{a3}. An alternative to this model choice is presented
in Section~\ref{experiments}.

The algorithm MOA, for minimizing $F$ over $\Omega$ based on the
regularized model, is defined as follows. Note that the third
condition in Assumption~\ref{a3} implies that $M(\bar x, x) + \sigma'
d(\bar x, x)^{p+1}$ is bounded below in $\Omega$ for all $\sigma' \geq
\sigma$.\\

\noindent
\textbf{Algorithm \ref{regula}.1 (MOA).}  Let $\sigma_{\min} > 0$,
$\alphabis > 0$, and $x_0 \in \Omega$ be given. Initialize $j
\leftarrow 0$.
\begin{description}
\item[\textbf{Step 1.}] Set $\ell \leftarrow 1$ and choose $\sigma_{j,
  1} \in [0, 1] $ such that $M(x_j, x) + \sigma_{j, 1} d(x_j,
  x)^{p+1}$ is bounded below in $\Omega$.

\item[\textbf{Step 2.}] Compute $x_{j,\ell} \in \Omega$ a solution of
  (\ref{surro}) with $\bar x = x_j$ and $\sigma = \sigma_{j, \ell}$.

\item[\textbf{Step 3.}] Test the condition
  \begin{equation} \label{armijoF}
    F(x_{j,\ell}) \leq F(x_j) - \alphabis d(x_j, x_{j,\ell})^{p+1}.
  \end{equation}
  If (\ref{armijoF}) does not hold, then set $\sigma_{j,\ell+1} =
  \max\{\sigma_{\min}, 2 \sigma_{j,\ell}\}$, $\ell \leftarrow \ell +
  1$, and go to Step~2.

\item[\textbf{Step 4.}] Set $x_{j+1} = x_{j,\ell}$, $\sigma_j =
  \sigma_{j,\ell}$, $j \leftarrow j+1$, and go to Step~1.\\
\end{description}

\noindent
\textbf{Remark~\ref{regula}.1.} It is worth noting that
Algorithm~\ref{regula}.1 may be seen as a generalization of the
projected gradient method for the convex constrained minimization of a
smooth function, in which a trial point $x_{j,\ell}$ of the form
$x_{j,\ell} = P_{\Omega}( x_j - \sigma_{j,\ell}^{-1} \nabla F(x_j))$
can be seen as the solution to subproblem~(\ref{surro}) with $M(\bar
x,x) := \nabla F(\bar x)^T (x-\bar x) + \sigma \| x - \bar x \|^2$,
$\bar x = x_j$, and $\sigma=\sigma_{j,\ell}$. Algorithm~\ref{regula}.1
has been defined without a stopping criterion. If $\eta > 0$ is a
small tolerance, condition
\[
d( x_{j+1},  x_j ) \leq \eta,
\]
tested at Step~4 right before going back to Step~1, would correspond
to the well-known stopping criterion
\begin{equation} \label{pg}
\|P_{\Omega}( x_j - \sigma_{j}^{-1} \nabla F(x_j)) - x_j \| \leq \eta,
\end{equation}
associated with the norm of the so called scaled continuous projected
gradient (evaluated at $x_j$). Note that criterion~(\ref{pg}) is
satisfied by~$x_j$ and that, by~(\ref{armijoF}), $F(x_{j+1}) \leq
F(x_j) - \gamma \eta^{p+1} < F(x_j)$. \\

The theorem below shows that Algorithm \ref{regula}.1 is well-defined
and, additionally, it gives an evaluation complexity result for each
iteration.

\begin{teo} \label{regula2}
Assume that $M(\cdot,\cdot)$, $L \geq 0$, and $p + 1 > 0$ are such
that Assumption~\ref{a3} holds. Then, the $j$th iteration of
Algorithm~\ref{regula}.1 is well defined and finishes with the
fulfillment of (\ref{armijoF}) after at most $O(\log(L+\gamma))$
evaluations of $F$.
\end{teo}

\begin{pro}
By (\ref{lipscho}), and the fact that, by the definition of
$x_{j,\ell}$, $M(x_j, x_{j,\ell}) + \sigma_{j, \ell} d(x_j ,
x_{j,\ell})^{p+1} \leq M(x_j, x_j)$, we have that
\[
\begin{array}{rcl}
F(x_{j,\ell}) &\leq& M(x_j, x_{j,\ell}) + L d(x_j,  x_{j,\ell})^{p+1} \\[2mm]
& = & M(x_j, x_{j,\ell}) + \sigma_{j,\ell}  d(x_j,  x_{j,\ell})^{p+1}  -
\sigma_{j,\ell} d(x_{j,\ell}, x_j)^{p+1} + L  d(x_j,  x_{j,\ell})^{p+1} \\[2mm]
&\leq& M(x_j, x_j)  + (L- \sigma_{j,\ell})  d(x_j,  x_{j,\ell})^{p+1} \\[2mm]
&=&  F(x_j)  + (L-   \sigma_{j,\ell})  d(x_j,  x_{j,\ell})^{p+1} .
\end{array}
\]       
Therefore, (\ref{armijoF}) holds if $\sigma_{j,\ell} \geq L+\alphabis$
that, by construction, occurs in the worst case when $\ell \geq
\log_{2} \left( (L+\alphabis)/\sigma_{\min} \right) + 2 $. Since $F$
is evaluated only at points $x_{j,\ell}$ to test condition
(\ref{armijoF}), this completes the proof.
\end{pro}

\begin{teo} \label{regula3}
Assume that $F(x) \geq F_{\mathrm{low}}$ for all $x \in \Omega$,
$M(\cdot,\cdot)$, $L \geq 0$, and $p + 1 > 0$ are such that
Assumption~\ref{a3} holds and the sequence $\{x_j\}$ is generated by
Algorithm~\ref{regula}.1. Then, given $\eta > 0$, the number of
iterations such that
\begin{equation} \label{mayorque}
d(x_j, x_{j+1}) > \eta
\end{equation}
is bounded above by
\begin{equation} \label{bounabove}
\left\lfloor \frac{F(x_0) - F_{\mathrm{low}}}{\alphabis \eta^{p+1}} \right\rfloor.
\end{equation}
Moreover, the number of evaluations of $F$ is bounded above by
\[
\left\lfloor \frac{F(x_0) - F_{\mathrm{low}}}{\alphabis \eta^{p+1}} \right\rfloor
\times \left( \log_{2} \left( \frac{L+\alphabis}{\sigma_{\min}} \right) + 2 \right).
\]  
\end{teo}

\begin{pro}
By (\ref{armijoF}) we have that, for all $j = 0, 1, 2,\dots$,
\begin{equation} \label{esta}
F(x_{j+1}) \leq F(x_j) - \alphabis d(x_j,  x_{j+1})^{p+1}.   
\end{equation}
If, for an iteration~$j$, \eqref{mayorque} holds, then, by~\eqref{esta}, 
\begin{equation} \label{esa}
F(x_j) - F(x_{j+1}) \geq \alphabis \eta^{p+1},
\end{equation}
meaning that, at this iteration~$j$, there is a decrease of the
objective function of at least $\alphabis \eta^{p+1}$. Since, by
definition, $x_j \in \Omega$ for all $j$ and $F(x) \geq
F_{\mathrm{low}}$ for all $x \in \Omega$, then the amount of such
decreases is limited by $(F(x_0) - F_{\mathrm{low}}) / (\alphabis
\eta^{p+1})$.  Thus the number of iterations at which~\eqref{mayorque}
holds is limited by~\eqref{bounabove}. The second part of the proof
follows from Theorem~\ref{regula2}.
\end{pro}

\begin{coro}
Assume that $F(x) \geq F_{\mathrm{low}}$ for all $x \in \Omega$,
$M(\cdot,\cdot)$, $L \geq 0$, and $p + 1 > 0$ are such that
Assumption~\ref{a3} holds and the sequence $\{x_j\}$ is generated by
Algorithm~\ref{regula}.1. Then,
\begin{equation} \label{limxk}
\lim_{j \to \infty}  d(x_j,  x_{j+1}) = 0.
\end{equation}
\end{coro}

\begin{pro}
The corollary follows immediately from Theorem~\ref{regula3}.
\end{pro}\\

Let us review the results obtained up to
now. Algorithm~\ref{mainalgo}.1 is our main Inexact Restoration
algorithm, for which we proved that $h(y_k) \to 0$ and $d(x_k,
x_{k+1}) \to 0$. Algorithm~\ref{solvinsub}.1 describes the
implementation of Step~3 of Algorithm~\ref{mainalgo}.1. However, for
the implementation of Algorithm~\ref{solvinsub}.1, we need a Monotone
Optimization Algorithm MOA (Algorithm~\ref{regula}.1), which was given
in this section. MOA is an iterative algorithm based on model
regularization that aims to minimize an arbitrary function $F$.  We
proved that the distance between consecutive iterates of MOA tends to
zero.  The question that we wish to answer now is: In which sense MOA
finds an approximate minimizer of $F$?  The answer to this question
requires definitions of criticality that are given below.

\begin{defi}
We say that $x_* \in \Omega$ is a critical point of the problem of
minimizing $F$ over $\Omega$, related to model $M(\cdot,\cdot)$ and $p
+ 1 > 0$ (in short, $(M,p)$-\textit{critical}) if there exists $\sigma
\in [0, 2 L]$ such that $x_*$ is a local minimizer of $ M(x_*, x) +
\sigma d(x_*, x)^{p+1}$ subject to $x \in \Omega$.
\end{defi}

In the following lemma we prove that local minimizers of $F$ over
$\Omega$ are, in fact, $(M,p)$-critical.

\begin{lem} \label{regula1}
Suppose that $x_* \in \Omega$ is a local minimizer of $F(x)$ subject
to $x \in \Omega$. Assume, moreover, that $M(\cdot,\cdot)$, $L \geq
0$, and $p + 1 > 0$ are such that Assumption~\ref{a3} holds. Then, for
all $\sigma \geq L$, $x_*$ is a local minimizer of $M(x_*, x) + \sigma
d(x_*, x)^{p+1}$ subject to $x \in \Omega$.
\end{lem}

\begin{pro}
Let $\sigma \geq L$ and suppose that there exists a sequence $\{x_j\}$
contained in $\Omega$ that converges to $x_*$ and $M(x_*, x_j) +
\sigma d(x_*, x_j)^{p+1} < M(x_*, x_*) = F(x_*)$ for all $j$.  Then,
since $\sigma \geq L$, $M(x_*, x_j) + L d(x_*, x_j)^{p+1} < F(x_*)$
for all $j$.  Therefore, by (\ref{lipscho}), $F(x_j) < F(x_*)$ for all
$j$ and, thus, $x_*$ is not a local minimizer of $F(x)$.
\end{pro}\\

The inspiration for the definition of criticality comes from the
analysis of the finite dimensional case in which the metric space
$\Omega$ is $\R^n$, $F$ is as smooth as desired, and $M(x_*, x)$ is
the Taylor polynomial of order $p$ around $x_*$, If $x_*$ is a local
minimizer it is natural to believe that this point will also be a
minimizer of $M(x_*, x)$. At least, this is what happens when $n=1$,
independently of the value of $p$, and when $p \in \{1, 2\}$,
independently of $n$. However, this is not true if $p > 2$ and $n >
1$. For example, take $F:\R^2 \to \R$ given by
\[
F(z_1, z_2) = z_2^2 - z_1^3 z_2 + z_1^6.
\]
The origin is a global minimizer of $F$. However, its Taylor
polynomial of order $p \in \{4,5\}$ is
\[
M(z_1, z_2) = z_2^2 - z_1^3 z_2,  
\]
for which the origin is not a local minimizer. For a similar example
with $p=3$, take
\[
F(z_1, z_2) = z_2^2 - z_1^2 z_2 + z_1^4. 
\]
In fact, when $x_*$ is a local minimizer of $F$, the regularized
Taylor polynomial has a minimum at $x_*$, if the regularization
parameter is not smaller than the Lipschitz constant~$L$.

The definition of critical point for the minimization of $F$ is
insightful but is not sufficient for practical purposes. As in most
optimization problems, we need a definition of approximate critical
points dependent of a given tolerance parameter $\eta>0$. This
definition is given below. In short, we say a point is $\eta$-critical
if its distance to a local minimizer of a regularized model is smaller
than $\eta$.

\begin{defi}
We say that $\bar x \in \Omega$ is an $\eta$-critical point, with
respect to the minimization of $F$ over $\Omega$, related to model
$M(\cdot,\cdot)$ and $p + 1 > 0$ (in short,
$\eta$-$(M,p)$-\textit{critical}) if there exists $\sigma \in [0,2L]$
and $z \in \Omega$ such that $d(\bar x, z) \leq \eta$ and $z$ is a
local minimizer of $M(\bar x, x) + \sigma d(\bar x, x)^{p+1}$ subject
to $x \in \Omega$.
\end{defi}

\begin{teo} \label{regula4}
Assume that $F(x) \geq F_{\mathrm{low}}$ for all $x \in \Omega$,
$M(\cdot,\cdot)$, $L \geq 0$, and $p + 1 > 0$ are such that
Assumption~\ref{a3} holds and the sequence $\{x_j\}$ is generated by
Algorithm~\ref{regula}.1. Then, given $\eta > 0$, the number of
iterations such that $x_j$ is not $\eta$-$(M,p)$-critical is bounded
above by
\[
\left\lfloor \frac{F(x_0) - F_{\mathrm{low}}}{\alphabis \eta^{p+1}} \right\rfloor.  
\]
\end{teo}

\begin{pro}
This proof follows immediately from Theorem~\ref{regula3}.
\end{pro}

\section{Condensed theoretical results} \label{condensed}

The following theorem summarizes the theoretical results related to
the solution of problem (\ref{theproblem}) using the techniques
described in this paper.  In the hypothesis of this theorem, we
included the decision that MOA is executed until the obtention of a
local minimizer of the regularized model only when $h(y_{k+1})$ is
small (smaller than or equal to $\varepsilon_{\feas}$). In fact, we do
not need to solve subproblems with high precision at previous
iterations of the main algorithm.  Of course, in practice it is
sensible to solve subproblems using MOA with reasonable precision,
otherwise our main purpose of being close to a solution when $h(y)$ is
small (and evaluation is expensive) could fail to be satisfied. It is
interesting to observe that even basic Inexact Restoration results
share this characteristic. See, for example, Remark~1
in~\cite[p.336]{ir3} with respect to the algorithm for the
optimization phase in classical constrained optimization problems.

\begin{teo} \label{condensado}
Consider the application of Algorithm~\ref{mainalgo}.1, with Step~3
implemented according to Algorithm~\ref{solvinsub}.1 and Step~3 of
Algorithm~\ref{solvinsub}.1 implemented according to
Algorithm~\ref{regula}.1 (MOA), for finding an approximate solution to
problem~(\ref{theproblem}). Suppose that Assumptions~\ref{a1},
\ref{a2}, and~\ref{a3} hold. Let $\varepsilon_{\feas},
\varepsilon_{\opt}, \eta > 0$ be given. Assume, in addition, that, if
$h(y_{k+1}) \leq \varepsilon_{\feas}$, then MOA is executed until it
finds an $\eta$-$(M,p)$-critical point for the problem of minimizing
$f(x, y_{k+1}) + \alpha d(x, x_k)^\nu$ subject to $x \in \Omega$.
Then there exists an iteration index~$k$ not larger than
\[
\frac{c_{\feas}}{\varepsilon_{\feas}} +
\frac{c_{\opt}}{\varepsilon_{\opt}^{\nu}},
\]
where $c_{\feas}$ only depends on $x_0$, $y_0$, $r$, $\theta_0$,
$\beta$, $h_{\max}$, and $f_{\min}$ and $c_{\opt}$ only depends on
$\alpha$, $x_0$, $y_0$, $r$, $\theta_0$, $\beta$, $h_{\max}$, and
$f_{\min}$, such that
\begin{description}
\item[(i)] $h(y_{k+1}) \leq \varepsilon_{\feas}$ and $d(x_k, x_{k+1})
  \leq \varepsilon_{\opt}$, and
\item[(ii)] $x_{k+1}$ is $\eta$-$(M,p)$-critical for the problem of
  minimizing $f(x, y_{k+1}) + \alpha d(x, x_k)^\nu$.
\end{description}
\end{teo}

\begin{pro}
Item (i) in the thesis follows from Theorems~\ref{mainalgo}.1
and~\ref{mainalgo}.2; while (ii) comes from Theorem~\ref{regula}.3 due
to the hypothesis on the execution of MOA when $h(y_{k+1}) \leq
\varepsilon_{\feas}$.
\end{pro}\\

Regarding the values of which~$c_{\feas}$ and~$c_{\opt}$ depend on,
note that~$x_0$, $y_0$, $\theta_0$, $\nu$, $r$, $\alpha$, and~$\beta$
are (constant) parameters of Algorithm~\ref{mainalgo}.1;
while~$h_{\max}$ and~$f_{\min}$ are problem dependent constants
referred to in Assumption~\ref{a2}. Recall that Assumption~\ref{a1}
says that restoration is always possible, Assumption~\ref{a2} says
that $h$ is bounded above and $f$ is bounded below, and
Assumption~\ref{a3} describes characteristics of the model that is
used in the MOA algorithm. Moreover, note that, if $h(y_{k+1}) >
\varepsilon_{\feas}$, it is enough to assume that the point $x$
computed by MOA satisfies $M(\bar x, x) + \sigma d(\bar x, x)^{p+1}
\leq M(\bar x, \bar x)$, i.e.\ that a ``rough'' minimization is done.

\section{Numerical experiments}  \label{experiments}

In this section, the methodology introduced in this work is applied to
the problem of reproducing a controlled experiment that mimics the
failure of a dam. The section begins by tackling the choice of an
adequate surrogate model and describing particular cases of
Algorithm~\ref{mainalgo}.1, Algorithm~\ref{solvinsub}.1 (Optimization
Phase), and Algorithm~\ref{regula}.1 (MOA algorithm for addressing
problem~(\ref{minefe}) at Step~3 of Algorithm~\ref{solvinsub}.1).

The choice of an adequate surrogate model is highly
problem-dependent. Whereas in the smooth case Taylor-like models
should be interesting, this is not the case when differentiability is
out of question. A default approach in the context of this paper
consists of choosing as model for $f(x, y_{k+1}) + \alpha d(x_k,
x)^\nu$ (see (\ref{Efe})) the function
\begin{equation} \label{default}
M(x_k, x) = [f(x, y_{k'}) + \alpha d(x_k, x)^\nu] - f(x_k, y_{k'}) + f(x_k, y_{k+1})
\end{equation}
where $k' < k+1$. In this way, $M(x_k, x_k) = F(x_k)$ as required at
Assumption~\ref{a3} and solving~(\ref{surro}) is, very likely,
affordable, as the evaluation of $f(x, y_{k'})$ should be cheaper than
the evaluation of $f(x, y_{k+1})$ due to the requirement $k' < k+1 $.
The particular case of Algorithms~\ref{mainalgo}.1--\ref{regula}.1
that employs~(\ref{default}) as a surrogate model follows.\\

\noindent
\textbf{Algorithm~\ref{experiments}.1}. Let $x_0 \in \Omega$, $y_0 \in
Y$, $\theta_0 \in (0,1)$, $\nu > 0$, $r \in (0,1)$, $\alpha > 0$,
$\beta > 0$, $\varepsilon_{\feas} > 0$, $\varepsilon_{\opt} > 0$, and
$\eta>0$ be given. Set $k \leftarrow 0$.
\begin{description}
\item[\textbf{Step 1.}] \textit{Restoration phase}

Define $\yrest \in Y$ satisfying (\ref{erre}) and (\ref{beta}).

\item[\textbf{Step 2.}] \textit{Updating the penalty parameter}
  
Compute $\theta_{k+1}$ by means of (\ref{cinco}) and~(\ref{seis}).

\item[\textbf{Step 3.}] \textit{Optimization phase}
  \begin{description}
  \item[\textbf{Step 3.1.}] Set $y_{k+1} = y_k$. 
  \item[\textbf{Step 3.2.}] If $h(y_{k+1}) \leq \varepsilon_{\feas}$,
    then call Algorithm~\ref{regula}.1 with the surrogate model
    defined by~(\ref{default}) returning when an
    $\eta$-$(M,p)$-critical point~$x$ is found. On return, define
    $x_{k+1}=x$ and stop.
  \item[\textbf{Step 3.3.}] Execute one iteration of
    Algorithm~\ref{regula}.1 with the surrogate model defined
    by~(\ref{default}) to find $\xtrial$. On return, if (\ref{67}) is
    fulfilled, then define $x_{k+1}=\xtrial$ and go to Step~4.
  \item[\textbf{Step 3.4.}] Set $y_{k+1} = \yrest$.
  \item[\textbf{Step 3.5.}] If $h(y_{k+1}) \leq \varepsilon_{\feas}$,
    then call Algorithm~\ref{regula}.1 with the surrogate model
    defined by~(\ref{default}) returning when an
    $\eta$-$(M,p)$-critical point~$x$ is found. On return, define
    $x_{k+1}=x$ and stop.
  \item[\textbf{Step 3.6.}] Execute one iteration of
    Algorithm~\ref{regula}.1 with the surrogate model defined
    by~(\ref{default}) to find $x$. On return, define $x_{k+1}=x$.
  \end{description}
\item[\textbf{Step 4.}] Set $k \leftarrow k+1$ and go to Step~1.
\end{description}

\noindent
\textbf{Remark~\ref{experiments}.1} If Algorithm~\ref{regula}.1 fails
in any of the tasks required at Steps~3.2, 3.3, 3.5, or 3.6 in
Algorithm~\ref{experiments}.1, this reveals that Assumption~\ref{a3}
does not hold. In this case, stop with an appropriate failure
message.\\

\noindent
\textbf{Remark~\ref{experiments}.2} If $h(y_{k+1}) >
\varepsilon_{\feas}$, we execute just one iteration of
Algorithm~\ref{regula}.1; see Steps~3.3 and~3.6. Note that, by
Theorem~\ref{condensado}, from the theoretical point of view, it is
not necessary to run Algorithm~\ref{regula}.1 up to high precision in
this case; therefore, to execute only one iteration is an admissible
choice. In fact, we expect that, in many cases, this is a reasonable
effort that we may invest in the cases that the objective function has
been computed poorly in order to approximate the current point to a
solution.\\

In what follows we describe the application of
Algorithm~\ref{experiments}.1 to the problem of fitting a descriptive
model for a controlled physical experiment that aims to simulate and
record the failure of a dam in \cite{pirulli}. As it is well known,
hydro-geological hazardous natural phenomena like debris flows,
a\-va\-lanches and submerged landslides can cause significant damage
and loss of lives and properties.  A number of tragic incidents all
over the world is well documented and analyzed in scientific
literature but these events are extremely complex and still
challenging in terms of mathematical modelling and numerical
simulations. The key interaction that yields fast flow events is the
one between solid grains and interstitial fluid. We briefly review
here the results reported in \cite{pirulli} which are based on a small
scale experiment and two numerical simulations with different methods.

The experiment \cite{pirulli} is a standard small-scale model of
column collapse with saturated material that allows propagation in air
in order to achieve similarity with natural flow-like landslide.  A
glass fume cuboid with fixed length and width is closed at one end
while the second end has a movable vertical gate (see Figure~1 of
\cite{pirulli}).  The cuboid is filled with a saturated granular
mixture.  The column height is varied to achieve different aspect
ratios between height and length.  The solid grains exhibit uniform
granular distribution and density. The gate and opening mechanism is
designed in such a way that partial desaturation prior to the collapse
is avoided and the fast uplift of the gate triggers the propagation of
the saturated mixture in a way that resembles natural disasters. The
results of this physical experiment are documented by high resolution
photos and essentially demonstrate similar behaviour for all aspect
ratios. However, the time evolution of the normalized front position
is different depending on the aspect ratio.  Initially, when the gate
is lifted both grain and water start moving forward at the open end of
column.  Then the upper parts head toward the bottom with the failing
surface evolving in time.  At the end of the process the granular
front stops while the water filters through the solid phase.

Two conceptually different approaches are addressed in \cite{pirulli}
for numerical simulations of column collapse.  The first one, based on
a Discrete Element Method (DEM), assumes that the granular material is
represented by an assembly of particles interacting at contact
points. The key issue, both conceptually and computationally, in this
approach is the definition of a contact model between particles.  The
second approach is based on continuum models that require the
definition of a constitutive model of the material, which is still an
open question for landslides. A number of solvers for continuum models
is available but they differ greatly in mathematical description and
discretization of the continuum domain. In \cite{pirulli} two
numerical methods are tested: a mixed continuum-discrete DEM-LBM
(Lattice-Boltzmann) method and a two-phase double point Material Point
Method (MPM).  DEM-LBM offers more insight into the microscopic
structure of the granular material but it has significant limitations
due to computational costs. Given that the cost increases with the
total number of particles and the decrease of mesh size, DEM-LBM is
suitable for large grain or small scale experiments.  The MPM approach
is more applicable to real-scale problems as its efficiency is not
influenced by the grain size. However, this method is very sensitive
to the choice of the constitutive model and in the simulations
reported in \cite{pirulli} predicted a faster collapse with higher
velocities than the DEM-LBM method.

\begin{figure}[ht!]
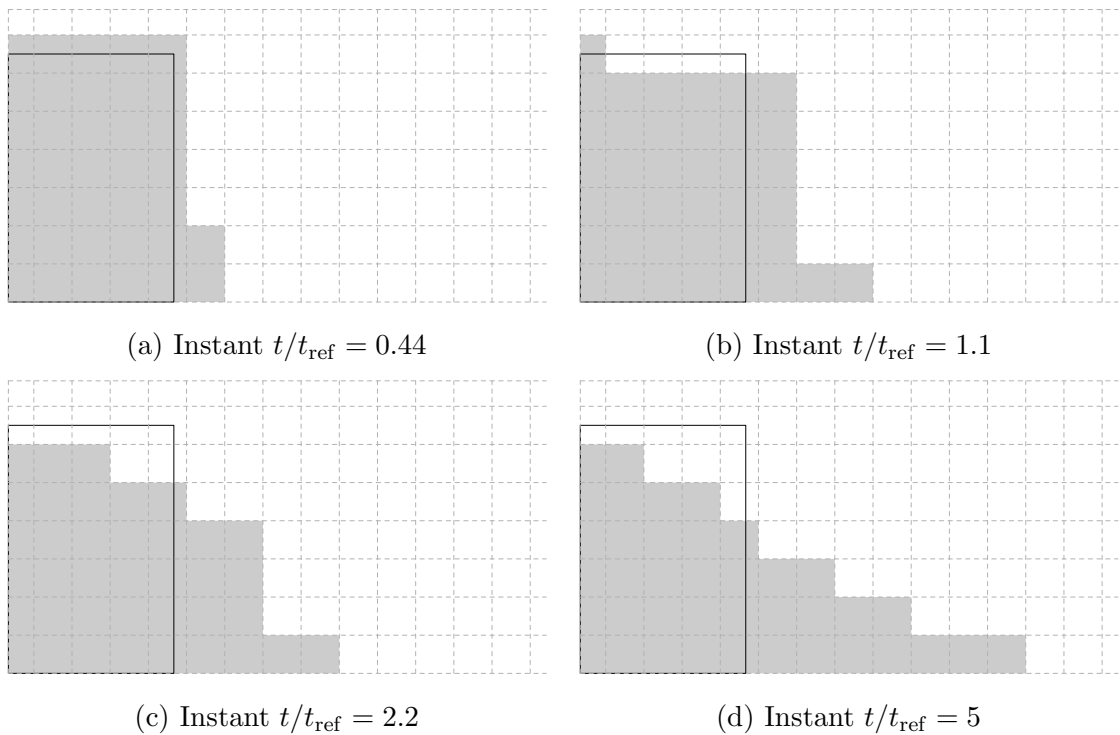

\begin{center}
\begin{tabular}{cc}
\includegraphics[scale=0.55]{bkm2020fig1-1.mps} &
\includegraphics[scale=0.55]{bkm2020fig1-2.mps}\\[2mm]
(a) Instant $t/t_{\mathrm{ref}}=0.44$ & 
(b) Instant $t/t_{\mathrm{ref}}=1.1$\\[2mm]
\includegraphics[scale=0.55]{bkm2020fig1-3.mps} &
\includegraphics[scale=0.55]{bkm2020fig1-4.mps}\\[2mm]
(c) Instant $t/t_{\mathrm{ref}}=2.2$ & 
(d) Instant $t/t_{\mathrm{ref}}=5$
\end{tabular}
\end{center}
\caption{Graphical representation of the four frames in~\cite[Fig.2,
    p.685]{pirulli}. Gray boxes represent boxes that are completely or
  partially filled with sand; while white boxes represent empty
  boxes.}
\label{fig1}
\end{figure}

In \cite{pirulli}, the considered granular material (sand) has a
uniform granulometric distribution with mean radius $0.125 \cm$ and
grain density $2625 \mathrm{km}/\mathrm{m}^3$. A rectangular box with
base $4\cm \times 5\cm$ and height~$7\cm$ is filled with sand
(approximately $8{,}000$ grains) and one of the $5\cm$-width sides is
removed. Four different frames corresponding to instants
$t/t_{\mathrm{ref}} = t_{\kappa}$ for $\kappa=1,2,3,4$ (with
$t_1=0.44$, $t_2=1.10$, $t_3=2.20$, and $t_4=5.0$) that illustrate the
dynamics of the system are presented. (See~\cite[Fig.2,
  p.685]{pirulli}.) Pictures have a squared paper on
background. Figures~\ref{fig1}(a--d) show, for each of the four frames
in~\cite[Fig.2, p.685]{pirulli}, the boxes completely or partially
filled with sand.

Let us now describe the approach we used to reproduce the above
described physical experiment. Denote $p_j=([p_j]_1,[p_j]_2)^T \in
\mathbb{R}^2$ for $j=1,\dots,n_p$, $p=(p_1^T,\dots,p_{n_p}^T)^T \in
\mathbb{R}^{2n_p}$, and
\begin{equation} \label{psidef}
\Psi_x(p) = x \left( \sum_{j=1}^{n_p} \sum_{i=j+1}^{n_p}
\max \left\{ 0, (2r)^2 - \| p_j - p_i \|_2^2 \right\}^2 \right) +
(1-x) \sum_{j=1}^{n_p} [p_j]_{2},
\end{equation}
where the ``radius'' $r>0$ and the ``weight'' $x \in (0,1) \subset
\mathbb{R}$ are given constant. In this formulation, $n_p$ represents
the number of balls, each $p_j \in \R^2$ ($j=1,\dots,n_p$) is the
center of a ball, the first term in~\eqref{psidef} (multiplied by $x$)
forces the balls' centers to be the ``correct distance'' apart; while
the second term in~\eqref{psidef} (multiplied by $1-x$) ``pushes all
the balls toward the ground''.

Consider the problem
\begin{equation} \label{prob}
\mbox{Minimize } \Psi_x(p_1,\dots,p_{n_p})
\mbox{ subject to } p \in D,
\end{equation}
where $D=\{ p \in \mathbb{R}^{2n_p} \; | \; p_j \geq 0 \mbox{ for }
j=1,\dots,n_p\}$. Independently of the value of $x \in (0,1)$, a
solution to problem~(\ref{prob}) consists in~$n_p$ points representing
the centers of $n_p$ non-overlapping identical balls with radius~$r$
``resting on the floor'' of the positive orthant of the
two-dimensional Cartesian space. Thus, the dam collapse is modeled not
just by changing $x$ but also by changing, simultaneously, the number
of iterations of an optimization methods applied to~\eqref{prob} as we
now describe.

\begin{figure}[ht!]
\begin{center}
\includegraphics[scale=0.55]{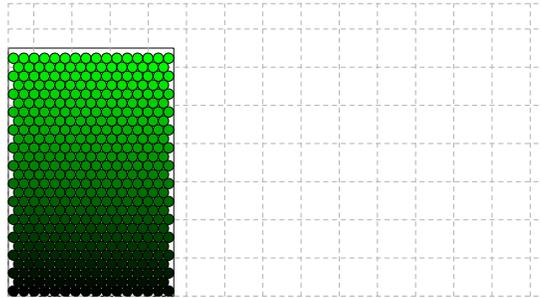}
\end{center}
\caption{Initial guess $p^0$ for the SPG method.}
\label{fig2}
\end{figure}

Assume that $p^0 \in D$ corresponds to the configuration represented
in Figure~\ref{fig2}. Given a maximum of iterations $y \in
\mathbb{N}_+$, let $p^0, p^1, p^2, \dots, p^{\bar y} \in D$ be
the~$\bar y \leq y$ iterates that result from the application of the
Spectral Projected Gradient (SPG) method~\cite{bmr,bmr2,bmr3,bmr5} to
problem~(\ref{prob}) starting from~$p^0$ and using as a stopping
criteria a maximum of~$y$ iterations or finding an iterate $p^{\bar
  y}$ such that
\[
\| P_D(p^{\bar y} - \nabla \Psi_x(p^{\bar y})) - p^{\bar y} \|_\infty
\leq \varepsilon_{\mathrm{opt}}^{\mathrm{spg}} := 10^{-8},
\]
where $P_D$ represents the projector operator onto $D$. We aim to
verify whether, by adjusting the weight~$x$ in~(\ref{psidef}) and the
maximum number of iterations~$y$ of the SPG method, it is possible to
construct a two-dimensional simulation of the dynamics of the
dam-failure physical experiment through the iterates $p^0, p^1, p^2,
\dots, p^{\bar y}$.

We consider solving problem~(\ref{prob}) with $r=0.125$ and $n_p=419$
using SPG with the initial guess $p^0$ depicted on
Figure~\ref{fig1}. We associate with the four frames depicted in
Figures~\ref{fig1}(a--d) binary matrices $M_1, M_2, M_3, M_4 \in
\{0,1\}^{8 \times 20}$ that represent whether there is sand in each
box of the frame or not. In an analogous way, we define $M(p) \in
\{0,1\}^{8 \times 20}$ as the matrix associated with the point
$p=(p_1^T,\dots,p_{n_p}^T)^T \in D \subset \mathbb{R}^{2n_p}$, that
indicates whether each box contains at least a point $p_j$ or
not. Given two $A=(a_{ij}), B=(b_{ij}) \in \{0,1\}^{8 \times 20}$, we
also define the fitness function
\[
\Pi(A,B) = \sum_{j=1}^{20} \sum_{i=1}^8 | a_{ij} - b_{ij} |.
\]
For a full sequence of iterates $p^0, p^1, p^2, \dots, p^{\bar y} \in
\Omega$ of the SPG method, that depends on the weight~$x$ and the
maximal number of iterations~$y$, we define
\[
f(x,y) = 1 - \frac{1}{640} \max_{c \geq 0} \left\{ \sum_{\kappa=1}^4
\Pi\left(M(p^{\lfloor c \, t_\kappa \rfloor}),M_\kappa\right) \right\},
\]
where, if $\lfloor c \, t_\kappa \rfloor > \bar y$ and, thus,
$p^{\lfloor c \, t_\kappa \rfloor}$ does not exist, then we consider
$\Pi\left(M(p^{\lfloor c \, t_\kappa \rfloor}),M_\kappa\right)=0$. The
value of $f$ can be computed by inspection of~$c$. Finally, we define
$h(y)=1/y$ and $\Omega=\{ x \in \R \;|\; 0 \leq x \leq 1 \}$.

\begin{figure}[ht!]
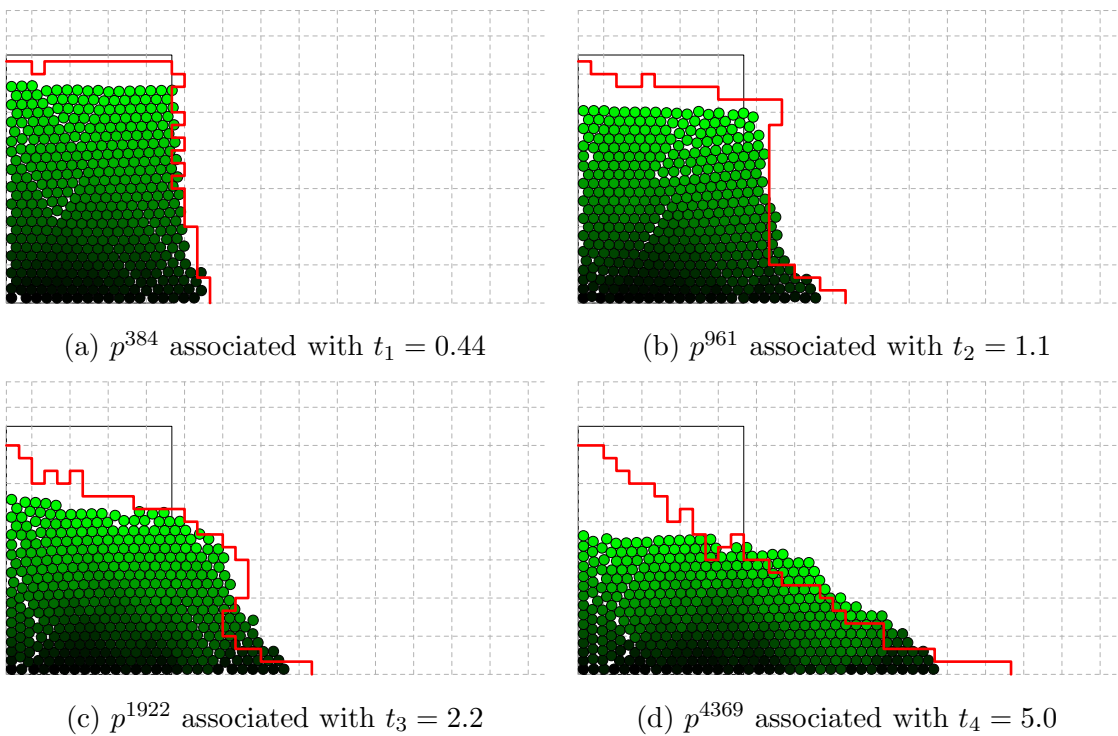

\begin{center}
\begin{tabular}{cc}
\includegraphics[scale=0.55]{bkm2020fig3-1.mps} &
\includegraphics[scale=0.55]{bkm2020fig3-2.mps}\\[2mm]
(a) $p^{384}$ associated with $t_1=0.44$ &
(b) $p^{961}$ associated with $t_2=1.1$ \\[2mm]
\includegraphics[scale=0.55]{bkm2020fig3-3.mps} &
\includegraphics[scale=0.55]{bkm2020fig3-4.mps}\\[2mm]
(c) $p^{1922}$ associated with $t_3=2.2$ &
(d) $p^{4369}$ associated with $t_4=5.0$
\end{tabular}
\end{center}
\caption{Graphical representation of $p^{384}$, $p^{961}$, $p^{1922}$,
  and $p^{4369}$ (corresponding to $c=873.2$) of a run of SPG with
  $x=0.999275$ and $y=12{,}800$. The red line is a rough
  representation of the physical experiment profile according to
  Fig. 2 of \cite{pirulli}.}
\label{fig3}
\end{figure}

We implemented Algorithm~\ref{experiments}.1 in Fortran. In the
numerical experiments, we considered, $\alpha=10^{-4}$, $\beta=100$,
$\theta_0=0.5$, $\nu = 2$, $r = 2$, $p=2$, $\varepsilon_{\feas} =
1/12{,}800$, $\varepsilon_{\opt} = 10^{-4}$, and $\eta=10^{-6}$. In
Algorithm~\ref{regula}.1, we considered $\sigma_{\min}=10^{-4}$ and
$\gamma=10^{-4}$. The optimality condition of Algorithm~\ref{regula}.1
defined by $\eta$ is that the objective function to which this
criterion is applied is not bigger, at the approximate optimizer $z$,
than the values at the feasible points of the form $z \pm \eta$. The
model $M(x_k,x)$ is given by~(\ref{default}) with $k' = \max\{0,
k-1\}$. For the approximate minimization of the model we employ
standard global one-dimensional search. As initial guess, we
considered $(x_0,y_0)=(0.5,100)$. All tests were conducted on a
computer with a 3.4 GHz Intel Core i5 processor and 8GB 1600 MHz DDR3
RAM memory, running macOS Mojave (version 10.14.6). Code was compiled
by the GFortran compiler of GCC (version 8.2.0) with the -O3
optimization directive enabled. Table~\ref{tab1} shows the results;
while Figure~\ref{fig3} illustrates de obtained solution. A rough
comparison with Figure~9 of \cite{pirulli} indicates that our results
are similar to the ones obtained by DEM-LBM and seem to be closer to
the experimental results than the ones reported for MPM. As in the
case of DEM-LBM and MPM, the discrepancies with respect to the
experimental results are due to the behavior of the top part of the
column.  The whole process was completed in less than 1 minute of CPU
time. The objective function equal to $1 - 618/640 \approx 0.03$ at
the final iterate means that at 97\% of the ``pixels'' of the measured
data matched the prediction of the model.

\begin{table}[ht!]
\begin{center}
\begin{tabular}{cccc}
\hline
$k$ & $x^k$ & $y^k$ & $f(x^k,y^k)$ \\
\hline
0 & 0.5      &      100 & $1-556/640$\\
1 & 0.9      &      200 & $1-561/640$\\
2 & 0.99     &      400 & $1-604/640$\\
3 & 0.99     &      800 & $1-604/640$\\
4 & 0.999    &  1{,}600 & $1-613/640$\\
5 & 0.99925  &  3{,}200 & $1-617/640$\\
6 & 0.999275 &  6{,}400 & $1-618/640$\\
7 & 0.999275 & 12{,}800 & $1-618/640$\\
\hline
\end{tabular}
\end{center}
\caption{Details of the application of Algorithm~\ref{experiments}.1.}
\label{tab1}
\end{table}

\section{Final remarks} \label{conclusions}

The application of Inexact Restoration to problems in which the
objective function is evaluated inexactly has proved to be successful
according to recent literature \cite{bkm2018,bkm2019,nkjmm}.
Real-life optimization problems use to be defined in domains more
general than finite-dimensional Euclidean spaces, so that
differentiability is out of question. The present paper aims to extend
the established theory to this more general context. For this purpose,
we found out that regularization-based algorithms are suitable for
solving subproblems that arise in the Optimization Phase of IR. The
regularization tools introduced in this paper are naturally inspired
in the Euclidean case, as this is the environment in which
regularization algorithms have been intensively developed in the last
few years. In particular, we needed to introduce new optimality
criteria and we emphasized the subtleties related with high-order
models for minimization problems. We do not claim that the approach
suggested in this paper for the Optimization Phase of IR is the
definite word on this subject. In fact, much research is necessary to
take advantage of the most reliable mathematical structures that are
associated with real-life problems. Topological vector spaces should
be considered in order to include problems in which the unknowns are
generalized functions (distributions) and the topology is given by
systems of semi-norms \cite{bourbaki}.

Many natural phenomena tend to converge to equilibrium states that are
characterized as minimizers of ``energy functions''. This is the
reason why physical laws sometimes inspire optimization
algorithms. For example, in \cite{zirilli} a method for solving
nonlinear equations is linked to a system of second-order ordinary
differential equations inspired in classical mechanics. In this paper,
we suggested a movement in the opposite direction: Simulating a
physical phenomenon by means of the behavior of an optimization
algorithm.  Of course this does not mean that physical events that are
usually modelled by means of complex systems of differential equations
may be always mimicked by the sequence of iterations of a simple
minimization algorithm. However, the possibility that a minimization
algorithm could give a first general approach to a complex phenomenon
in which the basic principle is energy minimization cannot be
excluded, especially when many physical parameters are unknown and, in
practice, need to be estimated using available data.

We are optimistic that the techniques described in this paper may help
understanding, describing and, at some point, predicting technological
disasters due to dam breaking. Many qualitative and partially
quantitative description of dam disasters exist in the scientific
literature.  For example, in \cite{picanco} we find a useful report
about the Brumadinho tailings dam disaster in Brazil in~2019. The
paper contains many references about this event, other dam disasters,
and the application of mathematical models to their analysis. See, for
example, \cite{stava}. We plan to apply the techniques introduced in
the present paper to the support of models of the Brumadinho event in
the context of the activity of CRIAB (acronym for ``Conflicts, Risks
and Impacts Associated with Dams'' in Portuguese), an
interdisciplinary research group created at the University of Campinas
for prediction and mitigation of the consequences of dam disasters.


\begin{thebibliography}{99}

\bibitem{ir7} R. Andreani, S. L. C. Castro, J. L. Chela,
  A. Friedlander, S. A. Santos, An inexact-restoration method for
  nonlinear bilevel programming problems, \textit{Computational
    Optimization and Applications} 43, pp. 307--328, 2009.
 
\bibitem{ir1} M. B. Aroux\'et, N. E. Echebest, and E. A. Pilotta,
  Inexact Restoration method for nonlinear optimization without
  derivatives, \textit{Journal of Computational and Applied
    Mathematics} 290, pp. 26--43, 2015.

\bibitem{ir13} N. Banihashemi and C. Y. Kaya, Inexact Restoration for
  Euler discretization of box-constrained optimal control problems,
  \textit{Journal of Optimization Theory and Applications} 156,
  pp. 726--760, 2013.

\bibitem{ir20} N. Banihashemi and C. Y. Kaya, Inexact restoration and
  adaptive mesh refinement for optimal control, \textit{Journal of
    Industrial and Management Optimization} 10, pp. 521--542, 2014.

\bibitem{inexact1} S. Bellavia, G. Gurioli, B. Morini, and
  Ph.L. Toint, Adaptive regularization algorithms with inexact
  evaluations for nonconvex optimization, \textit{SIAM Journal on
    Optimization} 29, pp. 2881--2915, 2019.
  
\bibitem{inexact2} S. Bellavia, G. Gurioli, B. Morini, and
  Ph.L. Toint, High-order Evaluation Complexity of a Stochastic
  Adaptive Regularization Algorithm for Nonconvex Optimization Using
  Inexact Function Evaluations and Randomly Perturbed Derivatives,
  arXiv preprint arXiv:2005.04639, 2020.
 
\bibitem{ir23} S. Bellavia, N. Kreji\'c, and B. Morini, Inexact
  restoration with subsampled trust‐region methods for finite‐sum
  minimization, \textit{Computational Optimization and Applications}
  76, pp. 701--736, 2020.

\bibitem{ir15} E. G. Birgin, L. F. Bueno, and J. M. Mart\'{\i}nez,
  Assessing the reliability of general-purpose Inexact Restoration
  methods, \textit{Journal of Computational and Applied Mathematics}
  282, pp. 1--16, 2015.
 
\bibitem{bkm2018} E. G. Birgin, N. Kreji\'c, and J. M. Mart\'{\i}nez,
  On the employment of Inexact Restoration for the minimization of
  functions whose evaluation is subject to errors, \textit{Mathematics
    of Computation} 87, pp. 1307--1326, 2018.
 
\bibitem{bkm2019} E. G. Birgin, N. Kreji\'c, J. M. Mart\'{\i}nez,
  Iteration and evaluation complexity on the minimization of functions
  whose computation is intrinsically inexact, \textit{Mathematics of
    Computation} 89, pp. 253--278, 2020.
 
\bibitem{ir24} E. G. Birgin, R. D. Lobato, and J. M. Mart\'{\i}nez,
  Constrained optimization with integer and continuous variables using
  inexact restoration and projected gradients, \textit{Bulletin of
    Computational Applied Mathematics} 4, pp. 55--70, 2016.
 
\bibitem{ir9} E. G. Birgin and J. M. Mart\'{\i}nez, Local convergence
  of an Inexact-Restoration method and numerical experiments,
  \textit{Journal of Optimization Theory and Applications} 127,
  pp. 229--247, 2005.

\bibitem{bmr} E. G. Birgin, J. M. Mart\'{\i}nez, and M. Raydan.
  Nonmonotone spectral projected gradient methods on convex sets.
  \textit{SIAM Journal on Optimization} 10, pp. 1196--1211, 2000.

\bibitem{bmr2} E. G. Birgin, J. M. Mart\'{\i}nez, and M. Raydan,
  Algorithm 813: SPG -- software for convex-constrained optimization,
  \textit{ACM Transactions on Mathematical Software} 27, pp. 340--349,
  2001.

\bibitem{bmr3} E. G. Birgin, J. M. Mart\'{\i}nez and M. Raydan,
  Inexact Spectral Projected Gradient methods on convex sets,
  \textit{IMA Journal of Numerical Analysis} 23, pp. 539--559, 2003.
  
\bibitem{bmr5} E. G. Birgin, J. M. Mart\'{\i}nez, and M. Raydan,
  Spectral Projected Gradient Methods: Review and Perspectives,
  \textit{Journal of Statistical Software} 60, issue 3, 2014.
 
\bibitem{bourbaki} N. Bourbaki, \textit{Elements of Mathematics ---
  Topological Vector Spaces, Chapters 1--5}, Springer, 2003.

\bibitem{bfms} L. F. Bueno, A. Friedlander, J. M. Mart\'{\i}nez, and
  F. N. C. Sobral, Inexact Restoration method for derivative-free
  optimization with smooth constraints, \textit{SIAM Journal on
    Optimization} 23, pp. 1189--1231, 2013.
 
\bibitem{ir16} L. F. Bueno, G. Haeser, and J. M. Mart\'{\i}nez, A
  flexible Inexact-Restoration method for constrained optimization,
  \textit{Journal of Optimization Theory and Applications} 165,
  pp. 188--208, 2015.

\bibitem{ir19} L. F. Bueno, G. Haeser, and J. M. Mart\'{\i}nez, An
  inexact restoration approach to optimization problems with
  multiobjective constraints under weighted-sum scalarization,
  \textit{Optimization Letters} 10, pp. 1315--1325, 2016.

\bibitem{ir22} L. F. Bueno and J. M. Mart\'{\i}nez, On the complexity
  of an Inexact Restoration method for constrained optimization,
  \textit{SIAM Journal on Optimization} 30, pp. 80--101, 2020.

\bibitem{inexact5} R. G. Carter, Numerical experience with a class of
  algorithms for nonlinear optimization using inexact function and
  gradient information, \textit{SIAM Journal on Scientific Computing}
  14, pp. 368--388, 1993.

\bibitem{pirulli} F. Ceccato, A. Leonardi, V. Girardi, P. Simonini,
  and M. Pirulli, Numerical and experimental investigation of
  saturated granular column collapse in air, \textit{Soils and
    Foundations} 60, pp. 683--696, 2020.

\bibitem{cgtbook} A. R. Conn, N. I. M. Gould, and Ph. L. Toint,
  \textit{Trust Region Methods}, MPS SIAM Series in Optimization,
  Philadelphia, 2000.
 
\bibitem{ir17} N. Echebest, M. L. Schuverdt, and R. P. Vignau, An
  inexact restoration derivative-free filter method for nonlinear
  programming, \textit{Computational and Applied Mathematics} 36,
  pp. 693--718, 2017.
 
\bibitem{ir11} D. Fern\'andez, E. A. Pilotta, and G. A. Torres, An
  inexact restoration strategy for the globalization of the sSQP
  method, \textit{Computational Optimization and Applications} 54,
  pp. 595--617, 2013.

\bibitem{ir18} P. S. Ferreira, E. W. Karas, M. Sachine,
  F. N. C. Sobral, Global convergence of a derivative-free inexact
  restoration filter algorithm for nonlinear programming,
  \textit{Optimization} 66, pp. 271--292, 2017.
 
\bibitem{ir3} A. Fischer and A. Friedlander, A new line search inexact
  restoration approach for nonlinear programming,
  \textit{Computational Optimization and Applications} 46,
  pp. 336--346, 2010.

\bibitem{ir25} J. B. Francisco, D. S. Gon\c{c}alves, F. S. V. Baz\'an,
  and L. L. T. Paredes, Non-monotone inexact restoration method for
  nonlinear programming, \textit{Computational Optimization and
    Applications} 76, pp. 867--888, 2020.

\bibitem{ir27} J. B. Francisco, D. S. Gon\c{c}alves, F. S. V. Baz\'an,
  and L. L. T. Paredes, Nonmonotone inexact restoration approach for
  minimization with orthogonality constraints, \textit{Numerical
    Algorithms} 86, pp. 1651--1684, 2021.
 
\bibitem{ir14} J. B. Francisco, J. M. Mart\'{\i}nez, L. Mart\'{\i}nez,
  and F. Pisnitchenko, Inexact restoration method for minimization
  problems arising in electronic structure calculations,
  \textit{Computational Optimization and Applications} 50,
  pp. 555--590, 2011.
 
\bibitem{ir12} M. A. Gomes-Ruggiero, J. M. Mart\'{\i}nez, and
  S. A. Santos, Spectral Projected Gradient method with Inexact
  Restoration for minimization with nonconvex constraints,
  \textit{SIAM Journal on Scientific Computing} 31, pp. 1628--1652,
  2009.

\bibitem{inexact7} S. Gratton, E. Simon, D. Titley-Peloquin,
  Ph. L. Toint, Minimizing convex quadratics with variable precision
  conjugate gradients, \textit{Numerical Linear Algebra with
    Applications} 28, e2337, 2021.
 
\bibitem{inexact4} S. Gratton, E. Simon, and Ph. L. Toint, An
  algorithm for the minimization of nonsmooth nonconvex functions
  using inexact evaluations and its worst-case complexity,
  \textit{Mathematical Programming} 187, 1--24, 2021.
 
\bibitem{inexact6} S. Gratton and Ph. L. Toint, A note on solving
  nonlinear optimization problems in variable precision,
  \textit{Computational Optimization and Applications} 76,
  pp. 917--933, 2020.
  
\bibitem{ir6} E. Karas, E. A. Pilotta, and A. Ribeiro, Numerical
  comparison of merit function with filter criterion in Inexact
  Restoration algorithms using hard-spheres problems,
  \textit{Computational Optimization and Applications} 44,
  pp. 427--441, 2009.

\bibitem{ir10} C. Y. Kaya, Inexact Restoration for Runge–Kutta
  discretization of optimal control problems, \textit{SIAM Journal on
    Numerical Analysis} 48, pp. 1492--1517, 2010.
 
\bibitem{ir8} C. Y. Kaya and J. M. Mart\'{\i}nez, Euler discretization
  and Inexact Restoration for optimal control, \textit{Journal of
    Optimization Theory and Applications} 134, pp. 191--206, 2007.

\bibitem{kltsiamreview} T. G. Kolda, R. M. Lewis, and V. Torczon,
  Optimization by direct search: new perspectives on some classical
  and modern methods, \textit{SIAM Review} 45, pp. 385--482, 2003.

\bibitem{klt} T. G. Kolda, R. M. Lewis, and V. Torczon, Stationarity
  results for generating set search for linearly constrained
  optimization, \textit{SIAM Journal on Optimization} 17,
  pp. 943--968, 2006.
 
\bibitem{inexact3} D. P. Kouri, M. Heinkenschloss, D. Ridzal, and
  B. G. van Bloemen Waanders, Inexact objective function evaluations
  in a trust-region algorithm for PDE-constrained optimization under
  uncertainty, \textit{SIAM Journal on Scientific Computing} 36,
  pp. A3011--A3029, 2014.

\bibitem{nkjmm} N. Kreji\'c and J. M. Mart\'{\i}nez, Inexact
  Restoration approach for minimization with inexact evaluation of the
  objective function, \textit{Mathematics of Computation} 85,
  pp. 1775--1791, 2016.
 
\bibitem{picanco} R. E. de Lima, J. L. Pican\c{c}o, A. F. da Silva,
  and F. A. Acordes, An anthropogenic flow type gravitational mass
  movement: the C\'orrego de Feij\~ao tailings dam disaster,
  Brumadinho, Brazil, \textit{Landslides} 17, pp. 2895--2906, 2020.

\bibitem{ir4} J. M. Mart\'{\i}nez, Inexact-restoration method with
  Lagrangian tangent decrease and new merit function for nonlinear
  programming, \textit{Journal of Optimization Theory and
    Applications} 111, 39--58, 2001.
 
\bibitem{ir2} J.M. Mart\'{\i}nez and E. A. Pilotta, Inexact
  restoration algorithms for constrained optimization, \textit{Journal
    of Optimization Theory and Applications} 104, pp. 135--163, 2000.

\bibitem{ir5} J. M. Mart\'{\i}nez and E. A. Pilotta, Inexact
  restoration methods for nonlinear programming: advances and
  perspectives, in L. Qi, K. Teo, X. Yang (eds.), \textit{Optimization
    and Control with Applications}, \textit{Applied Optimization} 96,
  Springer, New York, NY, 2005, pp. 271--291.

\bibitem{stava} M. Pirulli, M. Barbero, M. Martinelli, M. Marchelli,
  and C. Scavia, The failure of the Stava Valley tailings dam
  (Northern Italy): numerical analysis of the flow dynamics and
  rheological properties, \textit{Geoenvironmental Disasters} 4,
  article number 3, 2017.

\bibitem{ir21} C. E. P. Silva and M. T. T. Monteiro, A filter
  inexact-restoration method for nonlinear programming, \textit{TOP}
  16, pp. 126--146, 2008.

\bibitem{ir26} J. Walpen, P. A. Lotito, E. M. Mancinelli, and
  L. Parente, The demand adjustment problem via inexact restoration
  method, \textit{Computational and Applied Mathematics} 39, article
  number 204, 2020.

\bibitem{zirilli} F. Zirilli, The solution of nonlinear systems of
  equations by second order systems of O.D.E. and linearly implicit
  A-stable techniques, \textit{SIAM Journal on Numerical Analysis} 19,
  pp. 800--815, 1982.

\end{thebibliography}
\end{document}